\documentclass[12pt]{article}
\usepackage{amsmath}
\usepackage{amssymb}
\newcommand{\Z}{\mathbb Z}

\newcommand{\Q}{\mathbb Q}

\newcommand{\C}{\mathbb C}
\newcommand{\R}{\mathbb R}
\newcommand{\p}{\partial}
\newcommand{\z}{\zeta}

\newcommand{\la}{\lambda}
\newcommand{\La}{\Lambda}
\newcommand{\f}{\phi}
\newcommand{\LL}{\mathfrak{L}}
\newcommand{\bQ}{\overline{\Q}}

\begin{document}

\title{{\bf An introduction to Mathieu subspaces}}
\author{Arno van den Essen}
\maketitle

\section{Introduction}

The aim of this lectures is to give an introduction to the theory of Mathieu subspaces. We will not treat
all topics in their most general setting, but will restrict to certain classes of "nice" rings. As we  shall see even
for these rings there is a lot we don't know yet!

Mathieu subspaces grew out of several attempts to study the Jacobian Conjecture. One of these attempts was a result
due to Olivier Mathieu in 1995 ([Ma]). He stated the following conjecture and showed that his conjecture
implies the Jacobian Conjecture:

\medskip

\noindent{\bf Mathieu Conjecture.} {\em Let $G$ be a compact connected real Lie group with Haar measure $\sigma$.
Let $f$ be a complex valued $G$-finite function on $G$, such that
$\int_G f^m\,d\sigma=0$ for all positive $m$. Then for every $G$-finite function $g$ on $G$
also $\int_G gf^m\,d\sigma=0$, for all large $m$.}

\medskip

\noindent Here a function $f$ is called $G$-{\em finite} if the $\C$-vector space generated by the elements of the orbit $G\cdot f$ is finite dimensional.

\medskip

Then in 2003 another event happened: Michiel de Bondt and the author improved upon 
a classical result of Bass, Connell, Wright and Yagzhev ([BCW],[Y]). This result asserts that in order to
prove or disprove the Jacobian Conjecture it suffices to study so-called {\em cubic-homogeneous} polynomial maps, i.e.\@
maps of the form
$$x+H=(x_1+H_1,\cdots,x_n+H_n)$$
\noindent where the $H_i$ are either zero or homogeneous of degree three. Furthermore they showed that for
such a map $x+H$ the condition for its Jacobian determinant to be a non-zero constant is equivalent to the nilpotency
of the Jacobian matrix of $H$.

What de Bondt and the author showed in [BE]  is that we may {\em additionaly assume} that the Jacobian matrix
of $JH$ is {\em symmetric}, which by Poincar\'e's lemma implies that $H$ equals the gradient of some homogeneous polynomial  $P$ of degree four. In other words the Jacobian Conjecture is reduced to the study of
polynomial maps of the form
$$x+\nabla P$$
\noindent where $P$ is a homogeneous polynomial (of degree four).

When Wenhua Zhao heared about this result, he investigated what the Jacobian condition means for such a map
and also what it means for such a gradient map to be invertible. This led him to the following
conjecture, which he showed to be {\em equivalent} to the Jacobian Conjecture:

\medskip

\noindent{\bf Vanishing Conjecture.} {\em Let $\Delta=\sum_i\p_{x_i}^2$ be the Laplace operator and let $P\in\C[x]$
be homogeneous. If $\Delta^m(P^m)=0$ for all positive $m$, then $\Delta^m(P^{m+1})=0$ for all large $m$.}

\medskip

Zhao observed the resemblance with Mathieu's conjecture and in fact could made this resemblance even stronger by showing that the Vanishing Conjecture is {\em equivalent}
to the following version of it ([EZ]):

\medskip

\noindent{\bf Vanishing Conjecture.} {\em If $P\in\C[x]$ is homogeneous and satisfies $\Delta^m(P^m)=0$ for all
positive $m$, then for each $Q\in\C[x]$ also $\Delta^m(QP^m)=0$ for all large $m$.}

\medskip

At this point one can wonder: what makes the Laplace operator so special? The surprising answer is: it is not so special
at all. In fact it is not difficult to prove that if the Vanishing Conjecture holds for the Laplace operator, it also holds
for all quadratic homogeneous operators with constant coefficients. This observation led Zhao to the following
remarkable conjecture ([Z2]):

\medskip

\noindent{\bf Generalized Vanishing Conjecture (GVC(n)).} {\em Let $\Lambda$ be a differential operator with constant coefficients, i.e.\@ $\Lambda\in\C[\p_1,\cdots,\p_n]$. If $P\in\C[x]$ is such that $\Lambda^m(P^m)=0$ for
all positive $m$, then for each $Q\in\C[x]$ also $\Lambda^m(QP^m)=0$ for all large $m$.}

\medskip

Since, as we observed above, the Jacobian Conjecture is equivalent to the special case that $\Lambda$ is the
Laplace operator, we can view the Generalized Vanishing Conjecture as a family of "Jacobian Conjectures": for
each $\Lambda\in\C[\p_1,\cdots,\p_n]$ we get a new "Jacobian Conjecture".

When I saw this conjecture for the first time I was very sceptical and wanted to find a counter example. However
untill today no such example has been found, but instead various cases of the GVC have been proved!

But even this conjecture is not the end of the story, it is just the beginning: when Wenhua Zhao tried to prove it
he observed that the main obstruction is the fact that the operator $\Lambda$ and the multiplication operator $P$
do not commute. Therefore he studied {\em principal symbol maps}, a useful tool in the study
of differential equations. This led him to an even {\em more general} conjecture, the so-called Image Conjecture.

To describe this conjecture consider the  polynomial ring in $2n$ variables  $\C[\z, x]:=\C[\z_1,\cdots,\z_n,x_1,\cdots,x_n]$.
On it we have the following set of $n$ commuting differential operators
$$D:=\{\p_{x_1}-\z_1,\cdots,\p_{x_n}-\z_n\}$$
\noindent Put

$$ImD:=\sum_{i=1}^n (\p_{x_i}-\z_i)\,\C[\z,X]$$

\noindent{\bf Image Conjecture (IC(n)).} {\em If $f\in\C[\z,x]$ is such that $f^m\in ImD$, for all positive $m$, then
for each $g\in\C[\z,x]$ also $gf^m\in ImD$ for all large $m$.}

\medskip

Then in [Z3] Zhao showed that for each $n$ the Image Conjecture $IC(n)$ implies the Generalized Vanishing Conjecture $GVC(n)$. Consequently we have the following chain of implications

$$IC(n), \mbox{  for all } n\geq 1$$
$$\Rightarrow GVC(n), \mbox{ for all } n\geq 1$$
$$\Rightarrow VC(n), \mbox{ for all } n\geq 1$$
$$\Leftrightarrow JC(n), \mbox{ for all } n\geq 1$$

\noindent where $JC(n)$ denotes the $n$-dimensional Jacobian Conjecture.

\medskip 

The great insight of Zhao was the introduction of a general framework in which all conjectures mentioned above
can be studied simultaneously: this is his theory of {\em Mathieu subspaces}. In the next sections we will give the main
definitions and discuss several properties of these spaces. Before doing that let's make two remarks. The first remark
concerns the generality of the definitions to be given
below: we will not  strive for the greatest generality. Instead we will focuss on examples.  For that reason
we will mostly be concerned with the study of Mathieu subspaces of commutative rings. However we will make some
remarks on possible generalizations to the non-commutative case.

The second remark concerns the name Mathieu subspaces, which was given by its inventor in [Z4]. In order to honor
Wenhua Zhao for his insight, I propose to change the name of these spaces into 
{\em Mathieu-Zhao spaces}, abbreviated by MZ-spaces. From now on I will use this new name!

\section{Mathieu-Zhao spaces}

Throughout this section we have the following notations: $k$ is a field, $A$ any $k$-algebra containing $1$ and $V\subseteq A$ a $k$-linear subspace of $A$.

\medskip

\noindent{\bf Definition.} {\em Let $a,b\in A$. We say that $ba^{\infty}\in V$ if there exists $N\geq 1$ such that $ba^m\in V$ for all $m\geq N$.}

\medskip

Now we can introduce the following crucial notions:

\medskip

\noindent{\bf Definition.} {\em i) $r(V)$, the radical of $V$, is the set of $a\in A$ such that $a^{\infty}\in V$.\\
\noindent ii) $sr(V)$, the left strong radical of $V$, is the set of $a\in A$ such that for each $b\in A$, $ba^{\infty}\in V$.}

\medskip

Since $1\in A$, it follows readily that $sr(V)\subseteq r(V)$. However in general the converse inclusion does not hold:
for example take $A=k[t]$, the univariate polynomial ring, and $V=k[t^2]$. Then clearly $t^2\in r(V)$ since all powers
of $t^2$ belong to $V$. However, if we take $b=t$, then we see that none of the elements $t.(t^2)^m=t^{2m+1}$
belong to $V$. So $sr(V)$ is sticktly contained in $r(V)$.

\medskip

\noindent{\bf Definition.} {\em $V$ is called a left Mathieu-Zhao space (left MZ-space) of $A$ if $sr(V)=r(V)$. If $A$
is commutative we just say $V$ is a MZ-space.}

\medskip

\noindent In other words, a $k$-subspace $V$ of a $k$-algebra $A$ is called a left MZ-space
if the following holds: if $a^{\infty}\in V$, then $ba^{\infty}\in V$ for each $b\in A$. 

From this description we see that a left ideal in $A$ is a left MZ-space. So left MZ-spaces
are generalizations of left ideals. Below we will give plenty of examples of MZ-spaces
which are not ideals. On the other hand MZ-spaces are much harder to understand than
ideals. For example it is easy to describe all ideals in the univariate poynomial ring $k[t]$, but
describing MZ-spaces of this ring is still far to complicated. For example the set $V$ of all $f\in \C[t]$
such $\int_0^1 f(t)dt=0$ is a MZ-space of $\C[t]$ with $r(M)=0$. Proving this fact is not at all straightforward. We will
give a proof of it in section three.

\subsection*{Some examples}

\noindent{\bf Example 1.}  Let $A=M_n(k)$. Then
$$V:=\{C\in M_n(k)\,|\,TrC=0\,\}$$
\noindent is a left MZ-space of $A$.
\medskip

\noindent{\em Proof.} Let $C\in r(V)$. Then $C^{\infty}\in V$, i.e\@ $Tr\,C^m=0$ for all large $m$. Hence,
as is well-known, it follows that $C$ is nilpotent. In particular $C^n=0$.
So $ BC^m=0$  for all $B\in M_n(k) \mbox{ and all} \,m\geq n$. In particular
$Tr\, BC^m=0$  for all $B\in M_n(k) \mbox{ and all } \,m\geq n$.
So $ BC^m\in V$ for all  $B\in M_n(k) \mbox{ and all }\, m\geq n$, whence
$ BC^{\infty}\in V, \mbox{ for all } B\in M_n(k)$, ie.\@ $C\in sr(V)$.

\medskip

\noindent{\bf Example 2.} Let $A=\C[t][t^{-1}]$ and for $\la\in\C$ let $D_{\la}:=\p_t+\la t^{-1}$. Put
$$V_{\la}:=D_{\la}(A)$$
\noindent Then $V_{\la}$ is a MZ-space of $A$ if and only if $\la\notin\Z$ or $\la=-1$.

\medskip

\noindent{\em Proof.} Observe that $D_{\la}(t^i)=(\la+i)t^{i-1}$. So, if $\la\notin\Z$, then $D_{\la}(A)=A$,
so $V_{\la}$ is a MZ-space of $A$. Now let $\la\in\Z$ and assume that $\la\neq -1$. Then $\la+1\neq 0$, hence $1=D_{\la}(t/(\la+1))\in V_{\la}$. Suppose that $V_{\la}$ is a MZ-space of $A$, then for every $b\in A$ also $b.1^m\in V_{\la}$
for all large $m$. Hence $b\in V_{\la}$ for all $b\in A$. However, this is a contradiction since $\la\in\Z$, say $\la=m$,
and clearly $t^{-m-1}\notin V_{\la}$, since $D_{\la}(t^{-m})=(\la+{-m})t^{-m-1}=0.t^{-m-1}$. Finally
let $\la=-1$. Then it follows from the equality $D_{-1}(t^i)=(i-1)t^{i-1}$ that $V_{-1}$ is the set
of Laurent polynomials whose constant term equals $0$, i.e.
$$V_{-1}=\{f=\sum_i f_i\, t^i\in\C[t^{-1},t]\,|\,f_0=0\}$$
\noindent It is proved in [DvK] that $r(V_{-1})=t\C[t]\cup t^{-1}\C[t^{-1}]$, which easily implies
that $V_{-1}$ is a MZ-space of $A$.

\medskip

\noindent{\bf Example 3.} Let $A=\C[t]$. Then $V=\{f(t)\in A\,|\,\int_0^{\infty} e^{-t} f(t)\,dt=0\}$
\noindent is a MZ-space of $A$. More precisely, $r(V)=0$.

\medskip

\noindent A proof of this fact will be given in the next section. The reader is challenged to give a new proof!
In fact this result is a special case of the following conjecture:

\medskip

\noindent{\bf Integral Conjecture.} {\em Let $B\subseteq\R$ be an open set and $\sigma$ a positive measure such that
$\int_B g(t)d\sigma$ is finite for every $g(t)\in\C[t]$. Let
$$V_{\sigma}=\{f(t)\in\C[t]\,|\,\int_Bf(t)\,d\sigma=0\}$$
\noindent Then $V_{\sigma}$ is a MZ-space of $\C[t]$.}

\subsection*{Two simple observations}

To conclude this section we make two observations which will be used below. The first one
concerns the definition of a MZ-space. According to our definition $V$ is a left MZ-space of $A$ if for all
$a\in A$ such that $a^{\infty}\in V$ and each $b\in A$ we have that $ba^{\infty}\in V$. However we can strengthen
the condition $a^{\infty}\in V$ to $a^m\in V$ for all positive $m$. More precisely:

\medskip

\noindent{\bf Proposition 1.} {\em If for all $a\in A$ such that $a^m\in V$ for all positive $m$ it follows that for each $b\in A$
$ba^{\infty}\in V$, then $V$ is a left MZ-space of $A$.}

\medskip

\noindent{\em Proof.} Let $a\in r(V)$. Then there exists some $N$ such that $a^m\in V$ for all $m\geq N$.
In particular all positive powers of $a^N$ belong to $V$. Now let $b\in A$. Then, according to our hypothesis, for each of the elements
$ba^i$, with $0\leq i<N$ there exists $r_i$ such that $ba^i(a^N)^m\in V$ for all $m\geq r_i$.
Taking $r$ to be the maximum of all $r_i$ we obtain that for all $0\leq i<N$ we have $ba^{Nm+i}\in V$
for all $m\geq r$. It follows that $ba^n\in V$ for all $n\geq rN$, i.e.\@ $ba^{\infty}\in V$. So $a\in sr(V)$.

\medskip

It follows from this result that the Image Conjecture can be reformulated by saying that $ImD$ is a MZ-space
of $\C[\z,x]$. Similarly Mathieu's conjecture can be reformulated as follows: let $A$ be the ring of all complex valued
$G$-finite functions on $G$ and $V$ the set of all $f\in A$ such that $\int_G f\,d\sigma=0$. Then $V$ is
a MZ-space of $A$.

\medskip

The second observation, which will be used in the  next section, is the following: let $I$ be a left ideal in $A$
such that $I\subseteq V$.

\medskip

\noindent{\bf Proposition 2.} {\em $V$ is a MZ-space of $A$ iff $V/I$ is a MZ-space of $V/I$.}

\medskip

\noindent{\em Proof.} Straightforward.

\section{Tools}

\subsection*{Idempotents}

In this section we will develop some tools which will be useful in the study of MZ-spaces. The first of these 
is the so-called idempotent theorem, which was formulated and proved by Zhao in [Z5]. 

To describe this result observe the following: suppose that $V$ is a MZ-space of $A$ and that $e$ is an idempotent of
$A$ wich is contained in $V$, i.e.\@ $e^2=e$. Then obviously $e^m=e$  is contained in $V$ for all $m\geq 1$.
Since $V$ is a MZ-space of $A$ it follows that for each $b\in A$ the element $be^m$ belongs to $V$ for
all large $m$. However $be^m=be$. So we get that for  each $b\in A$ also $be\in V$. In other words $Ae\subseteq V$.

The Idempotent theorem below asserts that the converse holds if all elements of $r(V)$ are algebraic over $k$.
(Recall that an element $a\in A$ is called {\em algebraic over $k$} if there exists some non-zero polynomial $p(t)\in k[t]$
such that $p(a)=0$.)

\medskip

\noindent{\bf Idempotent theorem.} {\em Let $V$ be such that all elements of $r(V)$ are algebraic over $k$.
Then $V$ is a MZ-space of $A$ if and only if $Ae\subseteq V$ for all idempotents $e$ which belong to $V$.}

\medskip

\noindent{\em Proof.} We have already proved $(\Rightarrow)$, even without the assumption that the elements of
$r(V)$ are algebraic over $k$. So let's prove $(\Leftarrow)$. Therefore let $a\in r(V)$. Then for some $N$ $a^m\in V$ for
all $m\geq N$. Let $0\neq p(t)\in k[t]$ with $p(a)=0$. Then $q(t):=t^Np(t)=t^nh(t)$ with $n\geq N$ and $h(0)\neq 0$.
Also $q(a)=0$. Let $u(t),v(t)\in k[t]$ with
$$u(t)t^n+v(t)h(t)=1\,\,\,\,\,\,\,\,(*)$$
\noindent Put $e(t)=t^nu(t)$ and $e=e(a)$. Multiply (*) by $t^n$ and substitute $t=a$. Using $q(a)=0$ we get 
$a^ne=a^n$. Multiply (*) by $e(t)$ and substitute $t=a$. We get $e^2=e$. Also $e\in V$, since $n\geq N$. Hence $Ae\subseteq V$ (by the hypothesis). Finally, if $b\in A$ and $m\geq n$, then $ba^m=ba^{m-n}(a^ne)\subseteq Ae\subseteq V$.

\subsection*{p-adic techniques}

To explain the second tool for  studying MZ-spaces we need to recall some well known facts from algebraic geometry and number theory.
First recall that a field $k$ is algebraically closed if every polynomial $f(t)\in k[t]$ of positive degree has a zero in $k$.
Examples are the fields $\C$ and $\bQ$, the set of all complex numbers which are algebraic over $\Q$ .

Suppose now that $k$ is an {\em algebraically closed field} and that $f_1(x_1,
\cdots,x_n),\\\cdots,f_N(x_1,\cdots,x_n)$
is a set of polynomials in $k[x_1,\cdots,x_n]$ such that $1\in \sum_i k[x_1,\cdots,x_n]f_i(x_1,\cdots,x_n)$. Then
obviously the $f_i$ have no common zero in $k^n$. Now the {\em Nullstellensatz} asserts that the converse is true as well, i.e.\@
if a set of polynomials $f_i$ have no common zero in $k^n$, then $1$ belongs to the ideal in $k[x_1,\cdots,x_n]$
generated by the $f_i$.

\medskip

The second fact concerns some properties of absolute values, in particular of absolute values on number fields. 

Let $k$ be  any field. An {\em absolute value} on $k$
is a map $|\,\,\,  |:k\rightarrow \R$ which has the following properties: $|x|\geq 0$, for all $x\in k$ and $|x|=0$ iff $x=0$.
Furthermore $|xy|=|x||y|$ and $|x+y|\leq |x|+|y|$ for all $x,y\in k$. If $|x+y|\leq max(|x|,|y|)$, then $|\,\,\,|$ is
called a {\em non-Archimedean} absolute value.
Obviously the map $|\,\,\,|:k\rightarrow\R$ such that $|x|=1$ for all $x\in k-\{0\}$ and $|0|=0$ is an absolute value on $k$,
called the {\em trivial absolute value}.

On $\Q$ we can describe all absolute values: first there is the Euclidean absolute value, which we know from calculus . But there are more: for each prime
number $p$ we get an absolute value by the formula
$$|\frac{a}{b}p^n|_p=p^{-n}$$
\noindent where $a$ and $b$ have no factor $p$. This absolute value is called the {\em $p$-adic} absolute value on $\Q$.
It is a theorem due to Ostrowski which asserts that every non-trivial absolute value on $\Q$ is equivalent to either
the Euclidean absolute value or one of the $p$-adic ones (two absolute values $|\,\,|_1$ and $|\,\,\,|_2$ are called
{\em equivalent} if there exists some $c>0$ such that $|x|_2=|x|_1^c$ for all $x\in k$).

Now let $K$ be a {\em number field}, i.e.\@ $\Q\subseteq K$ is a finite field extension. Let $|\,\,\,|_p$ be the $p$-adic
absolute value on $\Q$ defined above. Then it is not difficult to show that this absolute value can be extended (in
a finite number of ways) to an absolute value on $K$. For each $p$ we choose such an extension
and denote it again by $|\,\,\,|_p$. Furthermore one can show the following property: for each $a\in K$ there
exists a number $N$ such that $|a|_p\leq 1$, for all $p>N$.

\medskip

To illustrate how the results above can be used to show that certain spaces are MZ-spaces, we will give a proof
of the following statement mentioned in the previous section:

\medskip

\noindent{\bf Theorem.} {\em Let $V=\{f(t)\in\C[t]\,|\,\int_0^{1} f(t)\,dt=0\}$. Then $r(V)=0$.
So $V$ is a MZ-space of $\C[t]$.}

\medskip

\noindent{\em Proof.} i)  Let $L:C[t]\rightarrow\C$ be defined by
$$L(f)=\int_0^{1} f(t)\,dt$$
\noindent Suppose $r(V)\neq 0$. Then there exists $0\neq f\in\C[t]$ such that $L(f^m)=0$ for all large $m$, say for
$m\geq N$. Hence $d:=deg_t f(t)\geq 1$ and we may assume that $f=t^d+\sum_{i=0}^{d-1}c_t\,t^i$, for some $c_i\in\C$.\\
\noindent ii) Next we will show that there exists $0\neq F\in\bQ[t]$ with $L(F^m)=0$ for all $m\geq N$. Therefore
observe that
$$f^m=t^{md}+\sum_{i=0}^{md-1} g_i^{(m)}(c)\,t^i\,\,\,\,\,\,\,\,(1)$$
\noindent where $c=(c_0,\cdots,c_{d-1})$ and $g_i^{(m)}(C_0,\cdots,C_{d-1})\in\Z[C_0,\cdots,C_{d-1}]$, the polynomial ring in the variables $C_j$ over $\Z$. Applying the linear map $L$ to the equation in (1) we obtain that
$$0=L(f^m)=L(t^{md})+\sum_{i=0}^{md-1} g_i^{(m)}(c)L(t^i)\,, \mbox{ for all } m\geq N\,\,\,\,\,\,\,\,(2)$$
\noindent Now for each $m\geq N$ we define the polynomial
$$P_m(C):=L(t^{md})+\sum_{i=0}^{md-1} g_i^{(m)}(C)L(t^i)$$
\noindent Observe that $L(t^i)\in\bQ$ for all $i$ (in fact $L(t^i)=\int_0^{1}t^i\,dt=\frac{1}{i+1}$), so 
$P_m(C)\in\bQ[C_0,\cdots,C_{d-1}]$ for all $m\geq N$. Since by (2) $P_m(c)=0$ for all $m\geq N$ it follows
that $1\notin J$, the ideal in $\bQ[C_0,\cdots,C_{d-1}]$ generated by the $P_m(C)$, with $m\geq N$. So by the Nullstellensatz (applied to the field $\bQ$) it follows that there exists an element $a=(a_0,\cdots,a_{d-1})\in{\bQ}^d$
such that $P_m(a)=0$ for all $m\geq N$. But this means that
$$F:=t^d+\sum_{i=0}^{d-1}a_i\,t^i\in\bQ[t]$$
\noindent and that $L(F^m)=0$ for all $m\geq N$. So, replacing $f$ by $F$,  we may
assume that all coefficients $c_i$ of $f$ belong to $\bQ$.\\
\noindent iii) Using that $L(t^i)=\frac{1}{i+1}$ equation (2) becomes
$$0=\frac{1}{md+1}+\sum_{i=0}^{md-1} g_i^{(m)}(c)\frac{1}{i+1}\,, \mbox{ for all } m\geq N\,\,\,\,\,\,\,\,(3)$$
\noindent  By Dirichlet's prime number theorem there exist infinitely many prime numbers of the form $md+1$. 
Since all $c_j$ are algebraic over $\Q$, there exists a number field $K$ which contains all the $c_j$.
Now choose $m\geq N$ such that $p=md+1$ is prime and furthermore large enough to guarantee that
$|c_j|_p\leq 1$ for all $j$. Then $|g_i^{(m)}(c)|_p\leq 1$ for all $i$. Also observe that for all $0\leq i\leq md-1<p$
we have that $|\frac{1}{i+1}|_p=1$ (since $1\leq i+1<p$ and hence $p$ does not divide $i+1$).
Consequently 
$$|\sum_{i=0}^{md-1} g_i^{(m)}(c)\frac{1}{i+1}|_p\leq 1\,\,\,\,\,\,\,(4)$$
\noindent Finally look at equation (3) and observe that $|\frac{1}{md+1}|_p=|\frac{1}{p}|_p=p>1$,
which contradicts (4).
 
\medskip

To conclude this section we will give a proof of the statement given in Example 3 of the previous section.
More precisely we show:

\medskip

\noindent{\bf Theorem.}  {\em Let $A=\C[t]$. Then $V=\{f(t)\in A\,|\,\int_0^{\infty} e^{-t} f(t)\,dt=0\}$
\noindent is a MZ-space of $A$. More precisely, $r(V)=0$.}

\medskip

\noindent{\em Proof.}  Let $L:C[t]\rightarrow\C$ be definied by
$$L(f)=\int_0^{\infty} f(t)e^{-t}\,dt$$
\noindent Suppose $r(V)\neq 0$. Then there exists $0\neq f\in\C[t]$ such that $L(f^m)=0$ for all large $m$, say for
$m\geq N$. We may assume that $f$ is of the form
$$f=t^r+c_{r+1}t^{r+1}+\cdots+c_{r+d}t^{r+d}$$
\noindent with $d\geq 1$. So for each $m\geq N$ we can write
$$f^m=t^{rm}+\sum_{i=1}^{md}g_i^{(m)}(c)t^{rm+i}\,\,\,\,\,\,(1)$$
\noindent where $c=(c_{r+1},\cdots,c_{r+d})$. Then, applying $L$ to (1) and using that
$L(t^i)=\int_0^{\infty}t^i\,e^{-t}\,dt=i!$,  we obtain
$$0=(rm)!+\sum_{i=1}^{md}g_i^{(m)}(c)(rm+i)!\,\,\,\,\,\,(2)$$
\noindent Arguing as above we may assume that all $c_j$ belong to some number field $K$. Divide the relation
in (2) by $(rm)!$ and observe that for each $1\leq i\leq md$ there exists a positive integer $n_i$
such that $\frac{(rm+i)!}{(rm)!}=(rm+1)n_i$. So from equation (2) we deduce
$$-1=\sum_{i=1}^{md}g_i^{(m)}(c)(rm+1)n_i\,\,\,\,\,\,(3)$$
\noindent Now again choose $m$ large enough such that $rm+1=p$ is prime and $|g_i^{(m)}(c)|_p\leq 1$ for all $i$. Then
it follows that the $p$-adic absolute value of the right hand side of (3) is stricktly smaller then $1$ (since $|rm+1|_p=|p|_p=
\frac{1}{p}$). However $|-1|_p=1$, which gives a contradiction.

\section{MZ-spaces of the univariate polynomial ring with non-zero strong radical}

As already observed in the previous sections, describing all MZ-spaces of the univariate polynomial ring $k[t]$ is
a problem which, at the moment, is still to difficult. In particular we don't know yet how to describe the MZ-spaces
whose strong radical is zero. On the other hand, as we shall show now, we can completely characterize
all MZ-spaces of $k[t]$ whose strong radical is non-zero, in case the field $k$ is algebraically closed and
its characteristic is zero. This description is essentially  based on three facts, namely the idempotent theorem,
a classical  result on linear recurrence relations and the following description of the strong radical of a subspace
of $k[t]$.

Let $R$ be a commutative $k$-algebra and $V$ a $k$-subspace of $R$. We define $I_V$ as the {\em largest
ideal} of $R$ contained in $V$. In other words, $I_V$ is the sum of all ideals contained in $V$. Obviously this
ideal can be the zero ideal. However in this section we will mainly be interested in the case where $I_V$ is non-zero.
The following result is the starting point for describing MZ-spaces of $k[t]$ with non-zero strong radical:

\medskip

\noindent{\bf Proposition.} {\em Let $V$ be a $k$-subspace of $k[t]$. Then $sr(V)=r(I_V)$.}

\medskip

\noindent{\em Proof.} i) Let $a\in r(I_V)$. Then $a^N\in I_V$ for some $N\geq 1$. Hence, since $I_V$ is an ideal, it follows
that $k[t]a^N\subseteq I_V\subseteq V$. In particular, for each $b\in k[t]$ we have that $ba^m=(ba^{m-N})a^N\in V$ for all $m\geq N$. So $a\in sr(V)$, which implies that $r(I_V)\subseteq sr(V)$.\\
\noindent ii) Conversely, if $sr(V)=0$, then by i) $r(I_V)=sr(V)$ and we are done. So assume that $sr(V)\neq 0$.
Choose $0\neq f\in sr(V)$, say $d=deg\,f$. If $d=0$, then $1\in sr(V)$, which implies that $V=k[t]$. Hence $sr(V)=r(I_V)$. So assume that $d\geq 1$.  Then $V$ is a (free) $k[f]$-module with basis $1,t,\cdots,t^{d-1}$. Since $f\in sr(V)$ there
exists $N$ such that $t^if^m\in V$ for all $0\leq i\leq d-1$ and all $m\geq N$. So $t^ik[f]f^N\subseteq V$ for all $0\leq i\leq d-1$. Since $k[t]=\sum_{i=0}^{d-1}t^ik[f]$, it follows that $k[t]f^N\subseteq V$. Since $k[t]f^N$ is an ideal contained in $V$, we obtain that $f^N\in k[t]f^N\subseteq I_V$, which implies $f\in r(I_V)$.

\medskip

From now on we assume: $V$ is a $k$-subspace {\em  strictly contained in $k[t]$}, $k$ is an {\em algebraically closed field of characteristic zero} and $sr(V)\neq 0$.

Since $sr(V)=r(I_V)$  it follows that $I_V\neq 0$.
 We will write
$I$ instead of $I_V$. Since $V$ is strictly contained in $ k[t]$, $I=k[t]f$ for some  polynomial $f$ of positive degree.
Say
$$f=\prod_{\la\in\La}(t-\la)^{m(\la)}$$
\noindent where $\La$ denotes the set of different zeroes of $f$ in $k$ and $m(\la)$ denotes the multiplicity
of $\la$. From proposition 2 in section two we deduce that $V$ is a MZ-space of $k[t]$ if and only if
$M:=V/I$ is a MZ-space of the ring $A:=k[t]/I$. So we need to investigate MZ-spaces of $A$. 

The ring $A$ is finite dimensional over $k$, hence all its elements are algebraic over $k$. It then follows
from the idempotency theorem that $M$ is a MZ-space of $A$ if and only if for each  idempotent $e$ of $A$, which
is contained in $M$, the ideal $Ae$ is contained in $M$. Therefore we first need to study
the idempotents of $A$. To do so we use the Chinese remainder theorem. By this theorem
 the map
$$\f:A\rightarrow B:=\prod_{\la\in\La}k[t]/(t-\la)^{m(\la)}$$
\noindent given by $\f(g+I)=(g+(t-\la)^{m(\la)})_{\la\in\La}$ is an isomorphism. Consequently, if we understand
all idempotents of $B$, we also understand all idempotents of $A$. Now observe that $B$ is a direct product of the rings 
$B_{\la}:=k[t]/(t-\la)^{m(\la)}$ and that each $B_{\la}$ has only two idempotents, namely $0$ and $1$. It follows
that the elements $e_{\la}=(0,\cdots,0,1,0,\cdots,0)\in B$ (where the $1$ appears at  the component with index $\la$)
form an {\em orthogonal basis of idempotents of $B$}, i.e.\@ each $e_{\la}$ is an idempotent of $B$, 
$e_{\la}\cdot e_{\mu}=0$ for all $\la\neq\mu\in\La$ and each non-zero idempotent of $B$ is of the form $\sum_{\la\in\La^{'}}e_{\la}$,
for some non-empty subset $\La^{'}$ of $\La$. By the isomorphism $\f$ there exist $g_{\la}\in k[t]$, such that
$\f(g_{\la}+I)=e_{\la}$. Consequently the elements $g_{\la}+I$ form an orthogonal bases of idempotents of $A$.

\medskip

The following result will be applied to the ring $A$ and the elements $g_{\la}+I$ and plays a crucial role in the proof of the main theorem:

\medskip

\noindent{\bf Lemma 1.} {\em Let $R$ be a commutative ring which has an orthogonal
basis $E$ of idempotents. If $M$ is a MZ-space of $R$ which does not contain
any non-zero element of $E$, then $0$ is the only idempotent of $R$ in $M$.}

\medskip

\noindent{\em Proof.} Let $e$ be a non-zero idempotent of $R$. Then $e=\sum_i e_i$, a finite 
sum of non-zero elements $e_i\in E$. Now assume
that $e\in M$. Since $M$ is a MZ-space of $R$ there exists an $N$ such that
$e_1e^m\in M$ for all $m\geq N$. However $e_1e^m=e_1e=e_1^2=e_1$. So $e_1\in M$, which contradicts our
hypothesis. So $e\notin M$. Hence $0$ is the only idempotent in $M$.

\medskip

Since $I=k[t]f\subseteq V$, $d:=dim_k k[t]/V$ is finite. Observe that $d\geq 1$, since
by assumption $V\neq k[t]$. So there is a $k$-isomorphism $\psi:k[t]/V\rightarrow k^d$. Let $\pi$ be the canonical map from $k[t]$ to $k[t]/V$. Then $\LL:=\psi\circ\pi$ is a surjective
$k$-linear map from $k[t]$ to $k^d$ such that $V$=\,ker$\LL$. So $V$ can be written as the kernel of a $k$-linear map
$\LL$, which has the additional property that
 $k[t]f\subseteq$ ker$\LL$. For such spaces $V$ the following criterion decides if they are MZ-spaces of $k[t]$:

\medskip

\noindent{\bf Theorem 1.} {\em With the notations as above let $\LL:k[t]\rightarrow k^d$ be a $k$-linear map such 
$V$=ker $\LL$. Then $V$ is a MZ-space of $k[t]$ if and only if for each non-empty subset $\La^{'}$ of $\La$ we have that $\sum_{\la\in\La^{'}} \LL(g_{\la})\neq 0$}.

\medskip

\noindent{\em Proof.}  i) Assume that $V$ is a MZ-space of $k[t]$. Then $M:=V/I$ is a MZ-space of $A:=k[t]/I$. We saw above that the element $\overline{g_{\la}}:=g_{\la}+I\in A$ is an idempotent. Claim: it does not belong to $M$. Assume it does. Then 
$A\overline{g_{\la}}\subseteq M$. Since $I\subseteq V$ this implies that $k[t]g_{\la}\subseteq V$. However $I$ is the
largest ideal contained in $V$.  So in particular $g_{\la}\in I$, whence $\overline{g_{\la}}=0$. Since $g_{\lambda}+I$
corresponds to $e_{\lambda}$ it is non-zero by construction, contradiction. So indeed $\overline{g_{\la}}\notin M$.

Now let $\La^{'}$ be a non-empty subset of $\La$ and suppose that $\sum_{\la\in\La^{'}}\LL(g_{\la})=0$. Then
$g:=\sum_{\la\in\La^{'}} g_{\la}\in$ ker $\LL=V$.  So $g+I\in M$ is a non-zero idempotent (since it
corresponds to $\sum_{\la\in\La^{'}} e_{\la}$ in $B$). Then  lemma 1 gives a contradiction, since
we showed that none of the $g_{\la}+I$ belong to $M$.\\
\noindent ii) Conversely, assume that for each non-empty subset $\La^{'}$ of $\La$ we have that $\sum_{\la\in\La^{'}} \LL(g_{\la})\neq 0$. Since all non-zero idempotents of $A$ are of the form  $\sum_{\la\in\La^{'}} g_{\la}+I$ (as shown above), it follows that none of these elements  belong to $M=V/I$. Consequently $0$ is the only idempotent of $A$
contained in $M$. Since obviously $A0\subseteq M$, it follows from the idempotency theorem that $M$ is a MZ-space of $A$ and hence that $V$ is a MZ-space of $k[t]$. 

\section{The Main theorem}

In the main theorem of this section we will give a concrete description of all MZ-spaces of $k[t]$ whose strong radical is non-zero. This is done in the following way. First we remark that, as  we saw  in theorem 1 of the previous section,  all such spaces can be described as kernels of certain $k$-linear maps, having some additional properties. Next we give an  explicit description of  all $k$-linear maps $\LL:k[t]\rightarrow k^d$ which have a non-zero ideal $I$ in their kernel.
Subsequently we decribe when for such maps $I$ is the {\em largest}  ideal contained in their kernel. Finally, having this description, we will use theorem 1 to state and prove the main theorem.

\subsection*{Linear maps having a non-zero ideal in their kernel}

From now on we assume that $\LL:k[t]\rightarrow k^d$ is a $k$-linear map, such that $k[t]f\subseteq$ker $\LL$, where
$f=\prod_{\la\in\La} (t-\la)^{m(\la)}$ 
and $\La$ denotes a non-empty finite subset of $k$. Write $\LL=(L_1,\cdots,L_d)$. So each $L_i:k[t]\rightarrow k$ is a $k$-linear map
having $k[t]f$ in its kernel. Let $N$ be the degree of $f$ and write $f=t^N-c_1t^{N-1}+\cdots +-c_N$, with $c_j\in k$.

To describe the next lemma we need some more notation: for each $\la\in k$ we denote by $S_{\la}$ the substitution
map from $k[t]$ to $k$, sending $g(t)$ to $g(\la)$. Furthermore we denote by $D$ the $k$-linear map $t\p_t$ from $k[t]$ to itself and write $\p$ instead of $\p_t$. Finally we let $\La_*:=\La-\{0\}$.

\medskip

\noindent{\bf Lemma 2.} {\em Let $L:k[t]\rightarrow k$ be a $k$-linear map.  If  $k[t]f\subseteq ker\,L$, then for every $\la\in\La$ there exists a polynomial $P_{\la}\in k[T]$ of degree $< m(\la)$, such that
$$L=S_0\circ P_0(\p)+\sum_{\la\in\La_*}S_{\la}\circ P_{\la}(D).$$}
\medskip

\noindent{\em Proof.} Let $a_i=L(t^i)$ for each $i$. Then $L(t^if)=0$ implies that 
$$L(t^{N+i}-c_1t^{N+i-1}+\cdots+-c_Nt^i)=0$$
\noindent Hence
$$a_{N+i}=c_1a_{N+i-1}+\cdots +c_Na_i, \mbox{ for all } i\geq 0$$
\noindent It then follows from the theory of linear recurrence relations with constant coefficients that there exist
$b_i\in k$, for each $0\leq i<m(0)$ and $b_{\la,i}\in k$ for each $\la\in\La_*$ and $0\leq i<m(\la)$,  such that
$$a_n=\sum_{0\leq i<m(0)} b_i e_{i,n}+\sum_{\la\in\La_*}\sum_{0\leq i<m(\la)} b_{i,{\la}} n^i{\la}^n, \mbox{ for all } n\geq 0,$$
\noindent where $e_i$ denotes the (infinite) sequence $(0,\cdots,0,1,0,0,\cdots)$ i.e.\@
$e_{i,i}=1$ and $e_{i,j}=0$ if $j\neq i$. By using that $S_{\la}\circ D^i(t^n)=n^i{\la}^n$, $S_0\circ\p^i(t^n)=i!$
if $n=i$ and $0$ if $n\neq i$, and using that $L(t^n)=a_n$, we obtain
$$L=\sum_{0\leq i<m(0)} \frac{b_i}{i!} S_0\circ {\p}^i+\sum_{\la\in\La_*}\sum_{0\leq i<m(\la)} b_{i,{\la}}S_{\la}\circ D^i$$
\noindent which completes the proof.

\medskip

Conversely, it follows readily from the first part of the next lemma that each map $L$ given in lemma 2 contains $k[t]f$ in
its kernel.

\medskip

\noindent{\bf Lemma 3.} {\em Let $\la\in k^*$,  $i\geq 1$ and $u\in k[t]$. Then\\
\noindent i) $(S_{\la}\circ D^i)((t-\la)^jk[t])=(S_0\circ\p^i)(t^jk[t])=0$,
if $j>i$.\\
\noindent ii) $(S_{\la}\circ D^i)(t-\la)^iu\neq 0$, if $u(\la)\neq 0$ and $(S_0\circ\p^i)(t^iu)\neq 0$, if $u(0)\neq 0$.}

\medskip

Finally we are able to describe all MZ-spaces $V$ of $k[t]$ which have a non-zero strong radical.
As we have seen above the hypothesis that $sr(V)$ is non-zero is equivalent to the fact that the largest ideal contained in $V$ is of the form $I=k[t]f$, where $f$ has positive degree (we also assume that $V$ is strictly contained in $k[t]$).
We showed that $V$ can be written as
$V$=ker$\LL$, where $\LL=(L_1,\cdots,L_d):k[t]\rightarrow k^d$ is a $k$-linear map (for some $d\geq 1$).
Furthermore it follows from lemma 2 that each $L_i$ is of the form
$$L_i=S_0\circ {P_0}^{(i)}(\p)+\sum_{\la\in\La_*}S_{\la}\circ{P_{\la}}^{(i)}(D).\,\,\,\,\,(*)$$
\noindent Observe that this description does not take into account that $I$ is the largest ideal contained
in ker$\LL$, but only expresses that $I\subseteq ker\,L_i$ for each $i$. Therefore we still need to describe what  it means for $I$ to be the {\em  largest ideal} contained in ker$\LL$. This is the content  of the next lemma:

\medskip

\noindent{\bf Lemma 4.} {\em $I$ is the largest ideal contained in ker$\LL$  if for each $\la\in\La$
the maximum of the degrees of all $P_{\la}^{(i)}$ is equal to  $m(\la)-1$.}

\medskip

\noindent{\em Proof.} By the remark just above lemma 3 it follows that $I=k[t]f$ is contained in ker$\LL$, where $f=\prod_{\la\in\La}(t-\la)^{m(\la)}$. So it remains to see that $k[t]f$ is the largest ideal with this property. Let
$k[t]g$ be the largest ideal contained in ker$\LL$. Then, since $k[t]f\subseteq$ker$\LL$, we get $f\in k[t]g$. So $g$ divides $f$, 
say $g=\prod_{\la\in\La}(t-\la)^{e_{\la}}$, with $e_{\la}\leq m(\la)$ for all $\la$. We will show that $e_{\la}\geq m(\la)$ for all $\la$, which implies that $f=g$ and hence completes the proof.

So assume that $e_{\la}\leq m(\la)-1$ for some $\la$ and write $\La=\{\la_1,\cdots,\la_r\}$. We may assume that $\la=\la_1$. Since ker$\LL$ contains $k[t]g$, it also contains $k[t]h$, where
$$h=(t-\la_1)^{m(\la_1)-1}\prod_{i=2,\cdots,r}(t-\la_i)^{m(\la_i)}$$
\noindent By our hypothesis there exists an $i$ such that the degree of $P_{\la_1}^{(i)}$ is equal to $m(\la_1)-1$.
We claim that $L_i(h)\neq 0$. This clearly contradicts the fact that $\LL(h)=0$ and hence completes the proof.

To see that $L_i(h)\neq 0$ we distinguish between the cases $\la_1=0$ and $\la_1\neq 0$. If $\la_1=0$
it follows from lemma 3 that $P_{\la_1}^{(i)}(\p)h\neq 0$.
Also it follows from lemma 3 that $P_{\la_j}^{(i)}(D)h=0$ for all $2\leq j\leq r$. Hence $L_i(h)=P_{\la_1}(\p)h\neq 0$.
A similar argument gives that $L_i(h)=P_{\la_1}(D)h\neq 0$, if $\la_1\neq 0$.

\medskip

\noindent{\bf Main theorem.} {\em Let $V\subseteq k[t]$ be such that $sr(V)$ is non-zero. Then\\
\noindent i) $V$=ker$\LL$, where $\LL=(L_1,\cdots,L_d):k[t]\rightarrow k^d$
is a $k$-linear map, $d\geq 1$ and each $L_i$ is of the form (*), with the additional property that
for every $\la$ the maximum of the degrees of all $P_{\la}^{(i)}$ is equal to $m(\la)-1$.\\
\noindent ii)  $V$ is a MZ-space of $k[t]$
if and only if for every non-empty subset $\La^{'}$ of $\La$ there exists an $i$ such that
$$\sum_{\la\in\La^{'}}{P_{\la}}^{(i)}(0)\neq 0.$$}

\medskip

\noindent{\em Proof.} The first part of this theorem is shown above. So assume that $V$=ker$\LL$, with each $L_i$ as described in i).  By theorem 1 we get that $V$ is a MZ-space of $k[t]$ if and only if for each  non-empty subset $\La^{'}$ of $\La$ the element $\sum_{\la\in\La^{'}} g_{\la}$ does not belong to ker$\LL$, in other words that for some $i$,
 $\sum_{\la\in\La^{'}}L_i( g_{\la})\neq 0$. By lemma 5 below we have  that $L_i(g_{\la})=P_{\la}^{(i)}(0)$. This completes the proof.

\medskip

\noindent{\bf Lemma 5.} {\em Let $L$ be as in Lemma 2. Then $L(g_{\la})=P_{\la}(0)$.}

\medskip

\noindent{\em Proof.} By definition of $g_{\la}$ we have that $g_{\la}\equiv 0 \mod (t-\mu)^{m(\mu)}$ if $\la\neq\mu$
and  $g_{\la}\equiv 1\mod (t-\la)^{m(\la)}$. Since the operator $P_{\mu}(D)$ has degree smaller than $m(\mu)$, it
follows from lemma 3 that $S_{\mu}\circ P_{\mu}(D)(k[t](t-\mu)^{m(\mu)})=0$. In particular  $S_{\mu}\circ P_{\mu}(D)(g_{\la})=0$
if $\mu\neq\la$. So, for $\la\in\La_*$ we get that $L(g_{\la})=S_0\circ P_0(\p)(g_{\la})+S_{\la}\circ P_{\la}(D)(g_{\la})$.
Similarly, since $\la\neq 0$ and the degree of the operator $P_0(\p)$ has degree smaller than $m(0)$ we have that $S_0\circ P_0(\p)(g_{\la})=0$. So $L(g_{\la})=S_{\la}\circ P_{\la}(D)(g_{\la})$. Finally observe that from 
 $g_{\la}\equiv 1\mod (t-\la)^{m(\la)}$ and the fact that $D(1)=0$ we obtain that $L(g_{\la})=P_{\la}(0)$, which proves the case $\la\neq 0$. The proof for the case $\la=0$ is similar and is left to the reader.

\section{The Image Conjecture}

In the first lecture we encountered the Image Conjecture. We mentioned that it implies the Generalized Vanishing
Conjecture, which in turn implies the Jacobian Conjecture. All these conjectures were formulated
for polynomial rings with complex coefficients.

In this section we generalize the Image Conjecture to polynomial rings with coefficients in commutative $k$-algebras, where
$k$ is an arbitrary field. The aim of this section is to give a proof of this more general Image Conjecture
in case the field $k$ has {\em positive characteristic $p$}.

So let's start describing this general setting: first $k$ is any field and $A$ is an arbitrary commutative $k$-algebra.
The polynomial ring in $n$ variables $x=(x_1,\cdots,x_n)$ over $A$ will be denoted by $A[x]$.
 Elements of the ring
$A[x]$ will simply be called polynomials, without refering to $A$ or $x$.
Let $a_1,\ldots,a_n$ be elements of $A$ and denote by $D$ the set of
commuting differential operators
$$\partial_{x_1}-a_1,\ldots,\partial_{x_n}-a_n.$$
\noindent Put
$$ImD=\sum_{i=1}^n (\partial_{x_i}-a_i)A[x].$$

\medskip

\noindent {\bf Image Conjecture (IC(n,A)).} {\em Assume that $(a_1,\ldots,a_n)$
is a regular sequence in $A$. Then $ImD$ is a MZ-space of $A[x]$}.

\medskip

A sequence $(a_1,\ldots,a_n)$ is
called a {\em regular sequence in $A$} if $a_1$ is no zero-divisor in $A$,
for each $i\geq 1$ the element $a_{i+1}$ is no zero-divisor
in $A/(a_1,\ldots,a_i)$ and the ideal generated by all $a_i$ is not equal
to $A$.

To get a feeling for the difficulty of the problem, the reader is 
invited to find a proof for the one dimensional case. In fact,
in this dimension the conjecture has only been proved in case the
ideal $Aa_1$ is a radical ideal. If additionaly $A$ is a UFD, also the non-radical
case has been proved. The Jacobian Conjecture follows from the very
special case where $A$ is the polynomial ring
$\C[\zeta_1,\ldots,\zeta_n]$ and $a_i=\zeta_i$ for each $i$.

\medskip

From now on (except in the crucial lemma below) we assume: 
{\em $k$ is a field of characteristic $p>0$}.
 The proof of the Image Conjecture under this assumption is based on the following
result, whose proof will be sketched at the end of this section
(we refer to [EWZ1] for more details).

\medskip

\noindent {\bf Crucial lemma.} {\em Let $k$ be any field and $b$ a 
polynomial of degree $d$. Denote by $b_d$ its homogeneous 
component of degree $d$. If $b$ belongs to $ImD$, then all 
coefficients of $b_d$ belong to the ideal $I$ of $A$ generated 
by all $a_i$.}

\medskip

\noindent {\bf Corollary.} {\em Let $f$ be a sum of monomials $f_a x^a$.
If $f^p$ belongs to $ImD$ then each $f_a^p$ belongs to $I$.}

\medskip

\noindent {\em Proof}. Write $f$ as a sum of homogeneous components
$f_i$. Then $f^p$ is a sum of the $f_i^p$. It then follows from the
crucial lemma that all coefficients of $f_d^p$ belong to $I$,
where $d$ is the degree of $f$. So $f_d^p$ is a sum of monomials
of the form $c a_i x^{ap}$, with $|a|=d$. Since each such a
monomial is equal to $(\p_i-a_i)(-cx^{ap})$, which belongs to $ImD$,
it follows that $f_d^p$ belongs to $ImD$. Substracting this
polynomial from $f^p$ we obtain that 

$$f_0^p+\ldots+f_{d-1}^p \in Im D.$$

\noindent Then the result follows by induction on $d$.

\medskip

\noindent {\em Proof of the Image Conjecture} (in case char $k=p>0$). 

\medskip

\noindent i) Since $\p^p=0$ on $A[x]$, we get
that $a_i^p x^a=(\p_i-a_i)^p(-x^a)\in ImD$. So 
a polynomial certainly belongs to $ImD$ {\em if} all its coefficients
belong to the ideal $J$ generated by all $a_i^p$.\\
\noindent ii) To show that $ImD$ is a MZ-space of $A[x]$ we need to show that
if $f$ is a polynomial such that $f^m\in ImD$ for all $m\geq 1$, then for each polynomial
$g$ there exists some $N$ such that $gf^m\in ImD$ for all $m\geq N$.  We will show the following
stronger statement:  if a polynomial $f$ is such that $f^p$ belongs to 
$ImD$ and $g$ is any polynomial, then $gf^m\in ImD$ for all $m\geq p^2$.  By i) it therefore suffices to show
that all coefficients of $gf^m$ belong to $J$,  if $m\geq p^2$. To understand these coefficients we
write $f$ as a sum of monomials $f_a x^a$. Since $f^p$ 
belongs to $ImD$, it follows from the corollary that each 
$f_a^p$ belongs to $I$ and hence each $f_a^{p^2}$ belongs to $J$. 
Since $f^{p^2}$ is a sum of the monomials $f_a^{p^2}x^{p^2a}$, 
it  readily follows that all coefficients of $gf^m$ belong
to $J$ if $m\geq p^2$. As observed above this concludes the proof.

\medskip

\noindent {\em Proof of crucial lemma (sketch).}

\medskip

\noindent We only sketch the case $n=2$. Let $b\in ImD$,
say 

$$b=(\p_1-a_1)p+(\p_2-a_2)q \,\,\,\,\,\,\,\,\,\,(*)$$

\noindent for some polynomials $p$ and $q$.
Let $d$ be the degree of $b$ and denote by $b_d$ the homogeneous
component of degree $d$. Now we assume for simplicity that the degrees
of both $p$ and $q$ are at most $d+2$. Then looking at the component
of degree $d+2$ in (*) we get

$$-a_1p_{d+2}-a_2q_{d+2}=0$$

\noindent where $p_i$ and $q_i$ denote the homogeneous components
of degree $i$ of $p$ and $q$ respectively. From the regularity hypothesis 
on the sequence $a_1,a_2$ it then follows that there exists a 
polynomial $g_{d+2}$, homogeneous of degree $d+2$, such that 

$$p_{d+2}=a_2g_{d+2} \mbox{ and } q_{d+2}=-a_1g_{d+2}.$$

\noindent Comparing the components of degree $d+1$ in the equation (*) we get

$$\p_1p_{d+2}+\p_2q_{d+2}-a_1p_{d+1}-a_2q_{d+1}=0.$$

\noindent Substituting the formulas for $p_{d+2}$ and $q_{d+2}$, found above, gives

$$a_1(p_{d+1}+\p_2g_{d+2})-a_2(q_{d+1}-\p_1g_{d+2})=0.$$

\noindent Again, from the regularity condition it then follows that there
exists some polynomial $g_{d+1}$, homogeneous of degree $d+1$, such that

$$p_{d+1}=-\p_2g_{d+2}+a_2g_{d+1} \mbox{ and } q_{d+1}=\p_1g_{d+2}+a_1g_{d+1}.$$

\noindent Comparing the components of degree $d$ in the equation (*) we get

$$b_d=\p_1p_{d+1}+\p_2q_{d+1}-a_1p_d-a_2q_d.$$

\noindent Finally, substituting the formulas for $p_{d+1}$ and $q_{d+1}$ in the 
last equality gives

$$b_d=-a_1(p_d-\p_2g_{d+1})-a_2(q_d-\p_1g_{d+1}).$$

\noindent which gives the desired result.

\end{document}